# *Coding theory applied to KU-algebras*

## **By**

*Samy M.Mostafa [a]*, Bayumy.A.Youssef[b], *Hussein A. Jad [c]*

## *Abstract*

The notion of a KU-valued function on a set is introduced and related properties are investigated. Codes generated by KU-valued functions are established. Moreover, we will provide an algorithm which allows us to find a KU-algebra starting from a given binary block code.

***Keywords***. KU-valued function, binary block code of KU-valued functions

--------------------------------------------------------------------------------------

## *1. Introduction*

BCK-algebras form an important class of logical algebras introduced by Iseki [5,6,7] and were extensively investigated by several researchers. The class of all *BCK*-algebras is a quasivariety. Iseki posed an interesting problem (solved by Wro*n*ski [13]) whether the class of *BCK*-algebras is a variety. In connection with this problem, Komori [9] introduced a notion of *BCC*-algebras and Dudek [1] redefined the notion of *BCC*-algebras by using a dual form of the ordinary definition in the sense of Komori. Dudek and Zhang [2] introduced a new notion of ideals in *BCC*-algebras and described connections between such ideals and congruences. C.Prabpayak and U.Leerawat ([11], [12]) introduced a new algebraic structure which is called KU-algebra. They gave the concept of homomorphisms of KU - algebras and investigated some related properties. These algebras form an important class of logical algebras and have many applications to various domains of mathematics, such as, group theory, functional analysis, fuzzy sets theory, probability theory, topology, etc. Coding theory is a very young mathematical topic. It started on the basis of transferring information from one place to another. For instance, suppose we are using electronic devices to transfer information (telephone, television, etc.). Here, information is converted into bits of 1's and 0's and sent through a channel, for example a cable or via satellite. Afterwards, the 1's and 0's are reconverted into information again. Due to technical problems, one can assume that while the bits are sent through the channel, there is a positive probability p that single bits are being changed. Thus the received bits could be wrong. The idea of coding theory is to give a method of how to convert the information into bits, such that there are no mistakes in the received information, or such that at least some of them are corrected. On this account, encoding and decoding algorithms are used to convert and reconvert these bits properly. One of the recent applications of BCK-algebras was given in the Coding theory [3,8 ,12]. In Coding Theory, a block code is an error-correcting code which encodes data in blocks. In the paper [8], the authors introduced the notion of BCK-valued functions and investigate several properties. Moreover,they established block-codes by using the notion of BCK-valued functions. they show that every finite BCK-algebra determines a block-code constructed a finite binary block-codes associated to a finite BCK-algebra. In [3,12] provided an algorithm which allows to find a BCK-algebra starting from a given binary block code.
In [12] the authors presented some new connections between BCK- algebras and binary block codes.
In this paper, we apply the code theory to KU- algebras and obtain some interesting results.

## 2. Preliminaries

Now, we will recall some known concepts related to KU-algebra from the literature which will be helpful in further study of this article.

***Definition 2.1*** [10,11] Algebra(X, ∗, 0) of type (2, 0) is said to be a KU -algebra, if it satisfies the following axioms:

$(ku_1)$ $0 * x = x$

$(ku_2)$ $x * y = 0 \Rightarrow (y * z) * (x * z) = 0, (z * x) * (z * y) = 0$

$(ku_3)$ $x * (y * z) = y * (x * z)$

$(ku_4)$ $(x * y) * [(y * z) * (x * z)] = 0$,

***Example 2.2*** Let X = {0, 1, 2, 3, 4} in which ∗ is defined by the following table

| ∗ | 0 | 1 | 2 | 3 | 4 |
|---|---|---|---|---|---|
| 0 | 0 | 1 | 2 | 3 | 4 |
| 1 | 0 | 0 | 2 | 3 | 4 |
| 2 | 0 | 1 | 0 | 3 | 3 |
| 3 | 0 | 0 | 2 | 0 | 2 |
| 4 | 0 | 0 | 0 | 0 | 0 |

It is easy to show that X is KU-algebra.

In a KU-algebra, the following identities are true : If we put in $(ku_4)$ y = x = 0 we get $(0 * 0) * [ (0 * z) * (0 * z) ] = 0$ and it follows that : $(Ku_5)$ $z * z = 0$, if we put y = 0 in $(ku_4)$, we get $(p_1)$ $z * (x * z) = 0$.

A subset S of KU-algebra X is called sub-algebra of X if $x * y \in S$, whenever x, y ∈ S.

A non empty subset A of a KU-algebra X is called a KU-ideal of X if it satisfies the following conditions:

$(I_1)$ 0 ∈ A,

$(I_2)$ $x * (y * z) \in A$, y ∈ A implies $x * z \in A$, for all x, y, z ∈ X.

***Lemma 2.3 [9]*** In a KU-algebra (X, ∗, 0), the following hold:

$x \leq y$ imply $y * z \leq x * z$.

***Lemma 2.4 [10]*** If X is KU-algebra then $y * [(y * x) * x] = 0$.

## 3. KU-valued functions

In what follows let A and X denote a nonempty set and a KU-algebra respectively, unless otherwise specified.

**Definition 3.1** A mapping $\tilde{A} : A \to X$ is called a KU-valued function (briefly, KU-function) on $A$.

**Definition 3.2** A cut function of $\tilde{A}$, for $q \in X$ is defined to be a mapping $\tilde{A}_q : A \to \{0,1\}$ such that $(\forall x \in A) \tilde{A}_q(x) = 1 \Leftrightarrow \tilde{A}(x) * q = 0$.

Obviously, $\tilde{A}_q$ is the characteristic function of the following subset of $A$, called a cut subset or a q-cut of $\tilde{A} : \tilde{A}_q(x) := \{ x \in A : \tilde{A}(x) * q = 0 \}$.

**Example 3.3** Let $A = \{x, y, z\}$ and let $X = \{0, a, b, c, d\}$ is a KU-algebra with the following Cayley table:

| * | 0 | a | b | c | d |
|---|---|---|---|---|---|
| 0 | 0 | a | b | c | d |
| a | 0 | 0 | b | b | a |
| b | 0 | a | 0 | a | d |
| c | 0 | 0 | 0 | 0 | a |
| d | 0 | 0 | b | b | 0 |

The function $\tilde{A} : A \to X$ given by $\tilde{A} = \begin{pmatrix} x & y & z \\ a & b & c \end{pmatrix}$ is a KU-function on $A$, and its cut subsets are

$A_0 = \Phi$ , $A_a = \{x\}$ , $A_b = \{y\}$, $A_c = A$ , $A_d = \{x\}$

**Proposition 3.4** Every KU-function $\tilde{A} : A \to X$ on $A$ is represented by the infimum of the set $\{ q \in X, \ A_q(x) = 1 \}$, that is $\forall x \in X \ : \tilde{A}(x) = \inf \{ q \in X, \tilde{A}_q(x) = 1 \}$.

**Proof.** For any $x \in A$. Let $\tilde{A}(x) = q \in X$, then $\tilde{A}(x) * q = 0$ and so $\tilde{A}_q(x) = 1 \ \forall q \in X$.

Assume that $\tilde{A}_r(x) = 1$ for $r \in X$, then $\tilde{A}(x) * r = 0 = q * r$, i.e $r \leq q$.

Since $q \in \{ r \in X, \tilde{A}_r(x) = 1\}$, for $x \in A$ , $r \in X$, it follows that $\tilde{A}(x) = q = \inf \{ r \in X, \tilde{A}_r(x) = 1 \}$.

This completes the proof.

***Proposition 3.5*** Let $\tilde{A}: A \to X$ be a KU-function on $A$. If $q * p = 0$ for all $p, q \in X$, we get $A_p \subseteq A_q$.

**Proof.** Let $p, q \in X$, be such that $q * p = 0$ and $x \in A_p$, then $\tilde{A}(x) * p = 0$

Using ($ku_1$) and ($ku_2$), we have

$$0 = \overbrace{(q * p) * (\tilde{A}(x) * P)}^{(KU_2)} = (\tilde{A}(x) * q),$$ and so $x \in A_q$. Therefore $A_p \subseteq A_q$.

This completes the proof.

***Proposition 3.6*** Let $\tilde{A}: A \to X$ be KU-function on $A$. Then
1- $(\forall x, y \in A)(\tilde{A}(x) \neq \tilde{A}(y) \Leftrightarrow A_{\tilde{A}(x)} \neq A_{\tilde{A}(y)}$

2- $(\forall q \in X)(x \in A)(\tilde{A}(x) * q = 0 \Leftrightarrow A_{\tilde{A}(x)} \subseteq A_q$

**Proof.** (1) The sufficiency is obvious. Assume that $A_{\tilde{A}(x)} \neq A_{\tilde{A}(y)}$ for all $x, y \in A$. Then

$A_{\tilde{A}(y)} * A_{\tilde{A}(x)} \neq 0$ or $A_{\tilde{A}(x)} * A_{\tilde{A}(y)} \neq 0$. Thus

$A_{\tilde{A}(x)} = \{ z \in A, \tilde{A}(z) * \tilde{A}(x) = 0 \} \neq \{ z \in A, \tilde{A}(z) * \tilde{A}(y) = 0 \} = A_{\tilde{A}(y)}$

(2) The necessity follows from Proposition 3.4. Let $q \in X$ and $x \in A$ be such that $A_{\tilde{A}(x)} \subseteq A_q$. If $\tilde{A}(x) * q \neq 0$ then $x \notin A_q$. Since $\tilde{A}(x) * \tilde{A}(x) = 0$, it follows that $x \in A_{\tilde{A}(x)}$, so that $A_{\tilde{A}(x)} \not\subset A_q$. This is a contradiction.

***Corollary 3.7*** Let $\tilde{A}: A \to X$ be KU-function on $A$. Then

$(\forall x, y \in A)(\tilde{A}(x) * \tilde{A}(y) = 0 \Leftrightarrow A_{\tilde{A}(y)} \subseteq A_{\tilde{A}(x)})$.

**Proof.** Straightforward.

For a KU-function $\tilde{A}: A \to X$, consider the following sets:
$A_x = \{ A_q : q \in X \}$, $\tilde{A}_x = \{ \tilde{A}_q : q \in X \}$.

***Proposition 3.8*** Let $\tilde{A}: A \to X$ be KU-function on $A$. Then

$(\forall Y \subseteq X)(\exists \inf Y$ in $X \Rightarrow A_{\inf(q: q \in Y)} = \bigcup \{ A_q : q \in Y \}$.

**Proof.** Let $(\forall Y \subseteq X)$ there exists $\inf Y$ in $X$ such that $x \in A_{\inf(q: q \in Y)}$. We have

$x \in A_{\inf(q: q \in Y)} \Leftrightarrow \tilde{A}(x) * \inf \{ q : q \in Y \} = 0 \Leftrightarrow (\forall r \in Y)(\tilde{A}(x) * r = 0) \Leftrightarrow (\forall r \in Y)(x \in A_r) \Leftrightarrow x \in \bigcup \{ A_q : q \in Y \}$. This completes the proof.

***Corollary 3.9*** Let $\tilde{A} : A \to X$ be KU-function on $A$, where $X$ is a bounded KU-algebra, then $\forall S \subseteq X$, $A_{\inf(q:\, q \in S)} = \bigcup \{A_q : q \in S\}$.

***Corollary 3.10*** Let $\tilde{A} : A \to X$ be KU-function on $A$, assume that for any $Y \subseteq X$, there exists a infimum of Y such that ( $\forall p, q \in Y$ ), we have $A_p \cup A_q \in A_X$.

The following example shows that the converse of the corollary 3.10 may not true in general.

**Example 3.11.** Let $A = \{x, y\}$ be a set and let $X = \{0, a, b, c, d\}$ be a KU-algebra with the following Cayley table:

| * | 0 | a | b | c | d |
|---|---|---|---|---|---|
| 0 | 0 | a | b | c | d |
| a | 0 | 0 | b | b | a |
| b | 0 | a | o | a | d |
| c | 0 | 0 | 0 | 0 | a |
| d | 0 | 0 | b | b | 0 |

The function $\tilde{A} : A \to X$ given by

$$\tilde{A} = \begin{pmatrix} x & y \\ a & b \end{pmatrix}$$ is a KU-function on $A$, then

|  | x | y |
|---|---|---|
|  | a | b |
| $\tilde{A}_0$ | 0 | 0 |
| $\tilde{A}_a$ | 1 | 0 |
| $\tilde{A}_b$ | 0 | 1 |
| $\tilde{A}_c$ | 1 | 1 |
| $\tilde{A}_d$ | 1 | 0 |

And its cut subsets are

$A_0 = \Phi$, $A_a = \{x\}$, $A_b = \{y\}$, $A_c = \{x, y\}$, $A_d = \{x\}$

Note that $A_a \cup A_b = \{x\} \cup \{y\} \notin A_X$, but $\inf\{a, b\}$ does exists in $X$.

***Proposition 3.12*** Let $\tilde{A}: A \to X$ be KU-function on $A$, then
$$\cap \{A_q | q \in X\} = A$$

**Proof.** Obviously, $\cap \{A_q | q \in X\} \subseteq A$. For every $x \in A$, let $\tilde{A}(x) = q \in X$. Then $x \in A_q$ and hence $x \in \cap \{A_q | q \in X\}$. Thus $A \subseteq \cap \{A_q | q \in X\}$. Therefore the result is valid.

***Proposition 3.13*** Let $\tilde{A}: A \to X$ be KU-function on $A$, then
$$(\forall x \in A)(\cup \{A_q | x \in A_q\} \in A_X)$$

**Proof.** Note that for any $x \in A, x \in A_q \Leftrightarrow \tilde{A}_q(x) = 1$,

From Proposition 3.7 we get the following

$\cup \{A_q | x \in A_q\} = \cup \{A_q | \tilde{A}_q(x) = 1\} = A_{\inf\{q | \tilde{A}_q(x)=1\}} \in A_q$. This completes the proof.

Let $\tilde{A}: A \to X$ be KU-function on $A$ and $\Theta$ be a binary operation on X defined by $\forall p, q \in X (p \Theta q \Leftrightarrow A_p = A_q)$. Then $\Theta$ is clearly an equivalence relation on X.

Let $\tilde{A}(A) = \{q \in X | \tilde{A}(x) = q \text{ for some } x \in A\}$ and for $q \in X$, $(q] = \{x \in X | x * q = 0\}$.

***Proposition 3.14*** For a KU-function $\tilde{A}: A \to X$ on $A$, we have
$$\forall p, q \in X (p \Theta q \Leftrightarrow (p] \cup \tilde{A}(A) = (q] \cup \tilde{A}(A)$$

**Proof.** We have $p \Theta q \Leftrightarrow A_p = A_q$
$$\Leftrightarrow (\forall x \in A)\{\tilde{A}(x) * p = 0 \Leftrightarrow \tilde{A}(x) * q = 0\}$$
$$\Leftrightarrow \{x \in A | \tilde{A}(x) \in (p]\} = \{x \in A | \tilde{A}(x) \in (q]\}$$
$$\Leftrightarrow (p] \cup \tilde{A}(A) = (q] \cup \tilde{A}(A).$$

This completes the proof.

***Example 3.15*** Let $X = \{a_n; n = 1, 2, 3, \ldots, 9\}$ and define a binary operation $*$ on $X$ as follows $(\forall a_i, a_j \in X)$ $(a_i * a_j = a_k)$, where $k = \dfrac{j}{(i,j)}$ and $(i, j)$ is the least common divisor of $i$ and $j$. Then $(X; *, a_i)$ is a KU-algebra. Its Cayley table is as follows:

| * | $a_1$ | $a_2$ | $a_3$ | $a_4$ | $a_5$ | $a_6$ | $a_7$ | $a_8$ | $a_9$ |
|---|---|---|---|---|---|---|---|---|---|
| $a_1$ | $a_1$ | $a_2$ | $a_3$ | $a_4$ | $a_5$ | $a_6$ | $a_7$ | $a_8$ | $a_9$ |
| $a_2$ | $a_1$ | $a_1$ | $a_3$ | $a_2$ | $a_5$ | $a_3$ | $a_7$ | $a_4$ | $a_9$ |
| $a_3$ | $a_1$ | $a_2$ | $a_1$ | $a_4$ | $a_5$ | $a_2$ | $a_7$ | $a_8$ | $a_3$ |
| $a_4$ | $a_1$ | $a_1$ | $a_3$ | $a_1$ | $a_5$ | $a_3$ | $a_7$ | $a_2$ | $a_9$ |
| $a_5$ | $a_1$ | $a_2$ | $a_3$ | $a_4$ | $a_1$ | $a_6$ | $a_7$ | $a_8$ | $a_9$ |
| $a_6$ | $a_1$ | $a_1$ | $a_1$ | $a_2$ | $a_5$ | $a_1$ | $a_7$ | $a_4$ | $a_3$ |
| $a_7$ | $a_1$ | $a_2$ | $a_3$ | $a_4$ | $a_5$ | $a_6$ | $a_1$ | $a_8$ | $a_9$ |
| $a_8$ | $a_1$ | $a_1$ | $a_3$ | $a_1$ | $a_5$ | $a_3$ | $a_7$ | $a_1$ | $a_9$ |
| $a_9$ | $a_1$ | $a_2$ | $a_1$ | $a_4$ | $a_5$ | $a_2$ | $a_7$ | $a_8$ | $a_1$ |

Let $A = \{a,b,c,d,e\}$ and $\tilde{A} : A \to X$ be a KU-function defined by

$$\tilde{A} = \begin{pmatrix} a & b & c & d & e \\ a_4 & a_6 & a_7 & a_1 & a_2 \end{pmatrix}. \text{ Then}$$

| * | a | b | c | d | e |
|---|---|---|---|---|---|
| | $a_4$ | $a_6$ | $a_7$ | $a_1$ | $a_2$ |
| $\tilde{A}_{a_1}$ | 0 | 0 | 0 | 1 | 0 |
| $\tilde{A}_{a_2}$ | 0 | 0 | 0 | 1 | 1 |
| $\tilde{A}_{a_3}$ | 0 | 0 | 0 | 1 | 0 |
| $\tilde{A}_{a_4}$ | 1 | 0 | 0 | 1 | 1 |
| $\tilde{A}_{a_5}$ | 0 | 0 | 0 | 1 | 0 |
| $\tilde{A}_{a_6}$ | 0 | 1 | 0 | 1 | 1 |
| $\tilde{A}_{a_7}$ | 0 | 0 | 1 | 1 | 0 |
| $\tilde{A}_{a_8}$ | 1 | 0 | 0 | 1 | 1 |
| $\tilde{A}_{a_9}$ | 0 | 0 | 0 | 1 | 0 |

and cut sets of Ã are as follows:

$\tilde{A}_{a_1} = \tilde{A}_{a_3} = \tilde{A}_{a_5} = \tilde{A}_{a_9} = \{d\}, \tilde{A}_{a_2} = \{d,e\}, \tilde{A}_{a_4} = \tilde{A}_{a_8} = \{a,d,e\}, \tilde{A}_{a_6} = \{b,d,e\}, \tilde{A}_{a_7} = \{c,d\}.$

## 4. Codes generates by KU-functions

Let $x/_\Theta = \{ y \in A ; \ x\Theta y \}$ ; for any $x \in A$, $x/_\Theta$ is called equivalence class containing $x$.

***Lemma 4.1*** Let $\tilde{A} : A \to X$ be a KU- function on $A$. For every $x \in A$, we have $\tilde{A}(x) = \inf \{x/_\Theta\}$, that is $\tilde{A}(x)$ the least element of the $\Theta$ to which it belongs.

**Proof.** Straightforward.

Let $A = \{1,2,3,\ldots,n\}$ and $X$ be a finite KU-algebra. Then every KU-function

$\tilde{A} : A \to X$ on $A$ determines a binary block code $V$ of length n in the following way: To every $x/_\Theta$, where $x \in A$, there corresponds a codeword $V_x = x_1 x_2 \ldots x_n$
Such that
$$x_i = x_j \Leftrightarrow \tilde{A}_x(i) = j \text{ for } i \in A \text{ and } j \in \{0,1\}.$$
Let $V_x = x_1 x_2 \ldots x_n$, $V_y = y_1 y_2 \ldots y_n$ be two code words belonging to a binary block-code $V$.

Define an order relation $\leq_c$ on the set of code words belonging to a binary block- code $V$ as follows: $V_x \leq_c V_y \Leftrightarrow x_i \leq y_i \ \text{for} \ i = 1,2,\ldots,n$ …… (4.1)

***Example 4.2*** Let $X = \{0, a, b, c\}$ be a KU-algebra with the following Cayley table:

| * | 0 | a | b | c |
|---|---|---|---|---|
| 0 | 0 | a | b | c |
| a | 0 | 0 | a | c |
| b | 0 | 0 | 0 | c |
| c | 0 | a | b | 0 |

Let $\tilde{A} : X \to X$ be a KU-function on $X$ given by
$\tilde{A} = \begin{pmatrix} 0 & a & b & c \\ 0 & a & b & c \end{pmatrix}$. Then

| $\tilde{A}_x$ | 0 | a | b | c |
|---|---|---|---|---|
| $\tilde{A}_0$ | 1 | 0 | 0 | 0 |
| $\tilde{A}_a$ | 1 | 1 | 0 | 0 |
| $\tilde{A}_b$ | 1 | 1 | 1 | 0 |
| $\tilde{A}_c$ | 1 | 0 | 0 | 1 |

$V = \{1000, 1100, 1110, 1001\}$. See Figure (1)

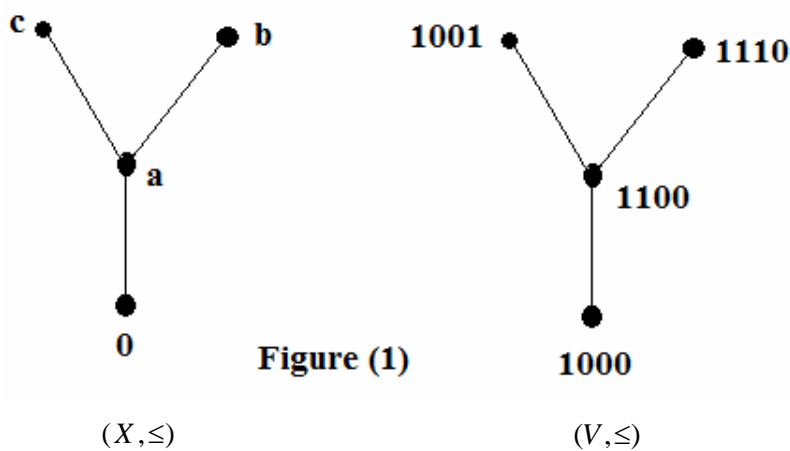

$(X, \leq)$          Figure (1)          $(V, \leq)$

Generally, we have the following theorem.

***Theorem 4.3*** Every finite KU-algebra X determines a block-code $V$ such that $(X, \leq)$ is isomorphic to $(X, \leq_c)$.

**Proof.** Let $X = \{a_i; i = 1,2,3,......,n\}$ be a finite KU-algebra in which $a_1$ is the least element and let $\tilde{A}: X \to X$ be identify KU-function on X. The decomposition of $\tilde{A}$ provides a family $\{\tilde{A}_q | q \in X\}$ which is the desired code under the order $V_x \leq_c V_y \Leftrightarrow x_i \leq y_i$ for $i = 1,2,....,n$

Let $f: X \to \{\tilde{A}_q; q \in X\}$ be a function defined by $f(q) = \tilde{A}_q$ for all $q \in X$. By lemma 4.1, every $\Theta$ class contains exactly one element. So, $f$ is one to one. Let $x, y \in X$ be such that $y * x = a_1$ i.e $x \leq y$. Then $A_x \subseteq A_y$ (by Proposition 3.5), which means that $\tilde{A}_x \subseteq \tilde{A}_y$. Therefore $f$ is an isomorphism. This completes the proof.

**Example 4.4**

Consider a KU-algebra $X = \{a_n; n = 1,2,3,......,9\}$ which is considered in example 3.15.
Let $\tilde{A}: X \to X$ be a KU-function on X given by

$$\tilde{A} = \begin{pmatrix} a_1 & a_2 & a_3 & a_4 & a_5 & a_6 & a_7 & a_8 & a_9 \\ a_1 & a_2 & a_3 & a_4 & a_5 & a_6 & a_7 & a_8 & a_9 \end{pmatrix}$$

Then

| * | $a_1$ | $a_2$ | $a_3$ | $a_4$ | $a_5$ | $a_6$ | $a_7$ | $a_8$ | $a_9$ |
|---|---|---|---|---|---|---|---|---|---|
| $\tilde{A}_{a_1}$ | 1 | 0 | 0 | 0 | 0 | 0 | 0 | 0 | 0 |
| $\tilde{A}_{a_2}$ | 1 | 1 | 0 | 0 | 0 | 0 | 0 | 0 | 0 |
| $\tilde{A}_{a_3}$ | 1 | 0 | 1 | 0 | 0 | 0 | 0 | 0 | 0 |
| $\tilde{A}_{a_4}$ | 1 | 1 | 0 | 1 | 0 | 0 | 0 | 0 | 0 |
| $\tilde{A}_{a_5}$ | 1 | 0 | 0 | 0 | 1 | 0 | 0 | 0 | 0 |
| $\tilde{A}_{a_6}$ | 1 | 1 | 1 | 0 | 0 | 1 | 0 | 0 | 0 |
| $\tilde{A}_{a_7}$ | 1 | 0 | 0 | 0 | 0 | 0 | 1 | 0 | 0 |
| $\tilde{A}_{a_8}$ | 1 | 1 | 0 | 1 | 0 | 0 | 0 | 1 | 0 |
| $\tilde{A}_{a_9}$ | 1 | 0 | 1 | 0 | 0 | 0 | 0 | 0 | 1 |

Thus
 V= {100000000, 110000000, 101000000, 110100000, 100010000, 11100100,100000100, 110100010, 101000001}. See Figure (2)

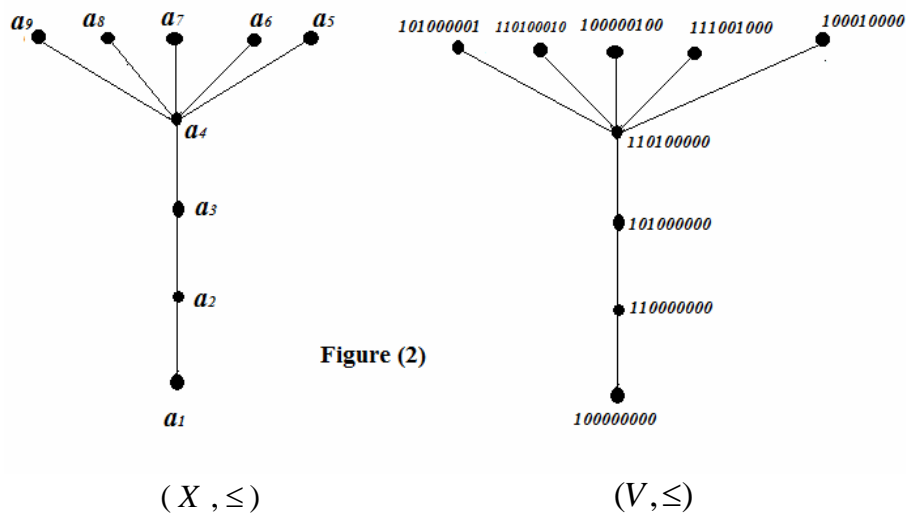

Figure (2)

$(X, \leq)$ $\qquad\qquad (V, \leq)$

Samy M. Mostafa   (samymostafa@yahoo.com)
Department of Mathematics, Faculty of Education, Ain Shams University, Roxy, Cairo, Egypt.

Bayumy.A.Youssef  (bbayumy@yahoo.com )
Informatics Research Institute,City for Scientific Research and Technological Applications, Borg ElArab, Alexandria, Egypt.

Hussein ali gad       (husseinaligad@yahoo.com )
Department of Mathematics and Computer Science .Informatics Research Institute(IRI)
City of Scientific Research and Technological applications, Alexandria, Egypt